\documentclass[10pt]{article}
\usepackage{amsfonts,amssymb,amsmath,amsthm,amsopn}
\theoremstyle{plain}
\newtheorem{teo}{Theorem}[section]

\newtheorem{lem}[teo]{Lemma}
\newtheorem{defi}[teo]{Definition}
\newtheorem{rem}[teo]{Remark}
\newcommand{\de}{\partial}
\renewcommand{\ge}{\varepsilon}

\begin{document}

\author{S. Secchi\thanks{Supported by M.U.R.S.T, Gruppo Nazionale 40 \% ``Variational Methods and Nonlinear Differential Equations''}\\
SISSA\\
via Beirut 2/4\\
34014 Trieste (Italy)}
\title{A Note on Closed Geodesics for a Class of non--compact Riemannian Manifolds}
\date{}
\maketitle

\section{Introduction}

\renewcommand{\thefootnote}{\fnsymbol{footnote}}

This paper is concerned with the existence of closed geodesics on a non--compact manifold $M$. There are very few papers on such a problem, see \cite{benci,tanaka,tho}. In particular, Tanaka deals with the manifod $M=\mathbb{R} \times S^{N}$, endowed with a metric $g(s,\xi)=g_0(\xi)+h(s,\xi)$, where $g_0$ is the standard product metric on $\mathbb{R} \times S^{N}$.
Under the assumption that $h(s,\xi)\to 0$ as $|s| \to \infty$, he proves the existence of a closed geodesic, found as a critical point of the energy functional
\begin{equation}
E(u)=\frac{1}{2} \int_0^1 g(u)[\dot u , \dot u] \, dt,
\end{equation}
defined on the loop space $\Lambda = \Lambda (M) = H^1 (S^1,M)$ \footnote[2]{We will identify $S^1$ with $[0,1] / \{ 0,1 \}$.}. The lack of compactness due to the unboundedness of $M$ is overcome by a suitable use of the concentration--compactness principle. To carry out the proof, the fact that $M$ has the specific form $M = \mathbb{R} \times S^{N}$ is fundamental, because this permits to compare $E$ with a {\em functional at infinity} whose behavior is explicitely known.

In the present paper, we consider a perturbed metric $g_\ge = g_0 + \ge h$, and extend Tanaka's result in two directions. First, we show the existence of at least $N$, in some cases $2N$, closed geodesics on $M=\mathbb{R} \times S^N$, see Theorem \ref{th:2.1}. Such a theorem can also be seen as an extension to cilyndrical domains of the result by Carminati \cite{carminati}. Next, we deal with the case in which $M=\mathbb{R} \times M_0$ for a general compact $N$--dimensional manifold $M_0$ \footnote[3]{By manifold we mean a smooth, connected manifold.}. The existence result we are able to prove requires that either $M_0$ possesses a non--degenerate closed geodesic, see Theorem \ref{th:3.1}, or that $\pi_1 (M_0) \ne \{0\}$ and the geodesics on $M_0$ are {\em isolated}, see Theorem \ref{th:4.1}.

\renewcommand{\thefootnote}{\arabic{footnote}}

The approach we use is different that Tanaka's one, and relies on a perturbation result discussed in \cite{ambad1} that leads to rather simple proofs. Roughly, the main advantages of using this abstract perturbation method are that
\begin{itemize}
\item[(i)] we can obtain sharper results, like the multeplicity ones;
\item[(ii)] we can deal with a general manifold like $M=\mathbb{R} \times M_0$, not only $M= \mathbb{R} \times S^N$, when the results -- for the reasons indicated before -- cannot be easily obtained by using Tanaka's approach.
\end{itemize}

The author wishes to thank professor Ambrosetti, for suggesting the problem, teaching him the perturbative technique, and much more.
\section{Spheres}

In this section we assume $M=\mathbb{R} \times S^N$, where $S^N=\{ \xi \in \mathbb{R}^{N+1} \colon | \xi | = 1 \}$ \footnote{Hereafter, we use the notation $\xi \bullet \eta = \sum_i \xi_i \eta_i$ for the scalar product in $\mathbb{R}^{N+1}$, and $|\xi|^2 = \xi \bullet \xi$.}.
For $s \in \mathbb{R}$, $r \in T_s \mathbb{R} \approx \mathbb{R}$, $\xi \in S^N$, $\eta \in  T_\xi S^N$, let
\begin{equation}
g_0(s,\xi)((r,\eta),(r,\eta))=|r|^2 + |\eta|^2
\end{equation}
be the standard product metric on $M=\mathbb{R} \times S^N$. We consider a perturbed metric
\begin{equation}
g_\ge (s,\xi)((r,\eta),(r,\eta))=|r|^2+|\eta|^2 + \ge h(s,\xi)((r,\eta),(r,\eta)),
\end{equation}
where $h(s,\xi)$ is a bilinear form, not necessarily positive definite.

Define the space of closed loops
\begin{equation} \label{Lambda}
\Lambda = \{ u = (r,x) \in H^1 (S^1,\mathbb{R}) \times H^1 (S^1 , S^N)\}
\end{equation}

Closed geodesics on $(M,g_\ge)$ are the critical points of $E_\ge \colon \Lambda \to \mathbb{R}$ given by
\begin{equation}\label{energia:perturbata}
E_\ge (u)=\frac{1}{2}\int_0^1 g_\ge (u)[\dot u , \dot u] \, dt.
\end{equation}
One has that
\begin{equation}\label{eq:vabene}
E_\ge (u)=E_\ge(r,x) = E_0(r,x)+\ge G(r,x),
\end{equation}
where
\begin{equation*}
E_0(r,x)=\frac{1}{2}\int_0^1 \left( |\dot r |^2 + | \dot x |^2 \right) \, dt
\end{equation*}
and
\begin{equation}\label{eq:G}
G(r,x)=\frac{1}{2} \int_0^1 h(r,x)[(\dot r,\dot x),(\dot r, \dot x)] \, dt.
\end{equation}
In particular, we split $E_0$ into two parts, namely
\begin{equation}
E_0(r,x)=L_0(r)+E_{M_0}(x),
\end{equation}
where
\[
L_0(r)=\frac{1}{2} \int_0^1 |\dot r|^2 \, dt, \quad E_{M_0}(x)=\frac{1}{2} \int_0^1 | \dot x|^2 \, dt.
\]

The form of $E_\ge$ suggests to apply the perturbative results of \cite{ambad1} that we recall below for the reader's convenience.

\begin{teo}[\cite{bertibolle,ambad1}]\label{th:berti:bolle}
Let $H$ be a real Hilbert space, $E_\ge \in C^2 (H)$ be of the form
\begin{equation}\label{eq:vabeneastratto}
E_\ge (u)=E_0(u)+\ge G(u),
\end{equation}
 where $G \in C^2 (E)$.. Suppose that there exists a finite dimensional manifold $Z$ such that
\begin{itemize}
\item[(AS1)] $E_0^\prime (z) = 0$ for all $z \in Z$;
\item[(AS2)] $E_0^{\prime \prime} (z)$ is a compact perturbation of the identity, for all $z \in Z$;
\item[(AS3)] $T_z Z= \ker E_0^{\prime \prime} (z)$ for all $z \in Z$.
\end{itemize}
There exist a positive number $\ge_0$ and a smooth function $w \colon Z \times (-\ge_0 , \ge_0) \to~H$ such that the critical points of
\begin{equation}\label{phi}
\Phi_\ge (z)=E_\ge (z+w(z,\ge)), \quad z \in Z,
\end{equation}
are critical points of $E_\ge$.
\end{teo}

Moreover, it is possible to show that
\begin{equation}\label{eq:gamma}
\Phi_\ge (z)=b+\ge \Gamma (z) + o(\ge),
\end{equation}
where $b=E_0 (z)$ and $\Gamma = G_{|Z}$.
From this ``first order'' expansion, one infers

\begin{teo}[\cite{ambad1}]
Let $H$ be a real Hilbert space, $E_\ge \in C^2 (H)$ be of the form \eqref{eq:vabeneastratto}. Suppose that
(AS1)--(AS3) hold.
Then any strict local extremum of $G_{|Z}$ gives rise to a critical point of $E_\ge$, for $|\ge|$ sufficiently small.
\end{teo}

\vspace{10pt}

In the present situation, the critical points of $E_0$ are nothing but the great circles of $S^N$, namely
\begin{equation}\label{zetapq}
z_{p,q}=p \cos 2 \pi t + q \sin 2 \pi t,
\end{equation}
where $p,q \in \mathbb{R}^{N+1}$, $p \bullet q = 0$, $|p|=|q|=1$.
Hence $E_0$ has a ``critical manifold'' given by
\begin{equation*}
Z=\{z(r,p,q)=(r,z_{p,q}(\cdot)) \mid r \in \mathbb{R}, \; z_{p,q} \mbox{\; as in \eqref{zetapq}} \}.
\end{equation*}

\begin{lem}
$Z$ satisfies (AS2)--(AS3).
\end{lem}

\begin{proof}
The first assertion is known, see for instance \cite{kli1}.

For the second statement, we closely follow \cite{carminati}.

For $z \in Z$, of the form $z(t)=(r,z_{p,q}(t))$, it turns out that
\[
 E_0^{\prime \prime} (z)[h,k]=\int_0^1 \left[ \dot h \bullet \dot k - |\dot z|^2 h \bullet k \right] \, dt
\]
for any $h,k \in T_z Z$.

Let $e_i \in \mathbb{R}^{N+1}$, $i=2, \dots,N+1$, be orthonormal vectors such that $\{ \frac{1}{2 \pi} {\dot z}_{p,q}, e_2, \dots e_{N+1}\}$ is a basis of $T_z Z$, and set
\begin{equation*}
e_i (t)=
\begin{cases}
\dot z_{u^1,u^2}(t) / 2 \pi &\text{ if } i=1 \\
e_i &\text{ if } i > 1,
\end{cases}
\end{equation*}

Then, for $h,k$ as before, we can write a ``Fourier--type'' expansion
\begin{equation} \label{acca}
h(t)=h_0 (t) \frac{d}{dt} + \sum_{i=1}^{N-1} h_i (t) e_i(t), \quad k(t)=k_0(t) \frac{d}{dt} + \sum_{i=1}^{N-1} k_i (t) e_i (t).
\end{equation}

Assume now that $h \in \ker E_0^{\prime \prime}(z_{p,q})$, i.e.
\[
\int_0^1 \dot h \bullet \dot k \, dt = \int_0^1 |\dot z|^2 h \bullet k \, dt \quad \forall k \in T_{z_{p,q}} Z .
\]

We plug \eqref{acca} into this relations, and we get the system

\begin{equation}
\begin{cases}
\overset{..}{h_1} = 0 \\
\overset{..}{h_j} + 4 \pi^2 h_j = 0 \quad j=2, \dots , N-1 \\
\overset{..}{h_0} = 0.
\end{cases}
\end{equation}

Recalling that $h_0$ and $h_1$ are periodic, we find
\begin{equation}
\begin{cases}
h_0 = \lambda_0, \quad h_1 = \lambda_1 \\
h_j = \lambda_j \cos 2 \pi t + \mu_j \sin 2 \pi t \quad j=2, \dots N-1 .
\end{cases}
\end{equation}

Therefore, $h \in T_z Z$. This shows that $\ker E_0^{\prime \prime}(z_{p,q}) \subset T_{z_{p,q}}Z$. Since the converse inclusion is is always true, the lemma follows.
\end{proof}

\begin{lem} \label{lem:due}
Suppose
\begin{itemize}
\item[(h1)] $h(r,\cdot) \to 0$ pointwise on $S^N$, as $|r| \to \infty$,
\end{itemize}
then
\[
\Phi_\ge \to b \equiv E_0(z).
\]
Recall that $\Phi_\ge$ was defined in \eqref{phi}.
\end{lem}

\begin{proof}
This is proved as in \cite{AM,bertibolle}. We just sketch the argument. The idea is to use the contraction mapping principle to characterize the function $w(\ge ,z)$ (see Theorem 1). Indeed, define
\begin{equation*}
H(\alpha , w,z_r,\ge )=
\begin{pmatrix}
E_\ge^\prime (z_r +w)-\alpha \dot z \\
w \bullet \dot z.
\end{pmatrix}
\end{equation*}
So $H=0$ if and only if $w \in (T_{z_r} Z)^\bot$ and $E_\ge^\prime (z_r+w) \in T_{z_r}Z$. Now,
\[
H(\alpha , w , z_r , \ge)=0 \Leftrightarrow H(0,0,z_r,0)+ \frac{\de H}{\de (\alpha ,w)} (0,0,z_r,0)[\alpha ,w] + R(\alpha ,w,z_r, \ge ) = 0,
\]
where $R(\alpha ,w,z_r, \ge ) = H(\alpha ,w,z_r,\ge)-\frac{\de H}{\de (\alpha ,w)} (0,0,z_r,0)[\alpha ,w]$.

Setting
\[
R_{z_r,\ge}(\alpha ,w)=- \left[ \frac{\de H}{\de (\alpha ,w)} (0,0,z_r,0) \right]^{-1}R(\alpha ,w,z_r, \ge),
\]
one finds that
\[
H(\alpha , w , z_r , \ge)=0 \Leftrightarrow (\alpha ,w) = R_{z_r,\ge}(\alpha ,w).
\]
By the Cauchy--Schwarz inequality, it turns out that $R_{zr,\ge}$ is a contraction mapping from some ball $B_{\rho (\ge)}$ into itself. If $|\ge|$ is sufficiently small, we have proved the existence of $(\alpha , w)$ uniformly for $z_r \in Z$. We want to study the asymptotic behavior of $w=w(\ge,z_r)$ as $|r| \to + \infty$. We denote by $R_\ge^0$ the functions $R_{z_r,\ge}$ corresponding to the unperturbed energy functional $E_0=E_{M_0}$. It is easy to see (\cite{bertibolle}, Lemma 3) that the function $w^0$ found with the same argument as before satisfies $\|w^0 (z_r)\| \to 0$ as $|r| \to + \infty$. Thus, by the continuous dependence of $w(\ge ,z_r)$ on $\ge$ and the characterization of $w(\ge ,z_r)$ and $w^0$ as fixed points of contractive mappings, we deduce as in \cite{AM}, proof of Lemma 3.2, that $\lim_{r \to \infty} w(\ge ,z_r)=0$. In conclusion, we have that $\lim_{|r| \to + \infty} \Phi_\ge (z_r + w(\ge ,z_r))=E_{M_0}(z_0)$.

\end{proof}

\begin{rem} \label{rem:3}
There is a natural action of the group $O(2)$ on the space $\Lambda$, given by
\begin{equation*}
\begin{align*}
\{ \pm 1 \} \times S^1 \times \Lambda &\longrightarrow \Lambda \\
(\pm 1 , \theta , u) & \mapsto u(\pm t + \theta ),
\end{align*}
\end{equation*}
under which the energy $E_\ge$ is invariant. Since this is an isometric action under which $Z$ is left unchanged, it easily follows that the function $w$ constructed in Theorem \ref{th:berti:bolle} is invariant, too.
\end{rem}

\begin{teo}\label{th:2.1}
Assume that the functions $h_{ij}=h_{ji}$'s are smooth, bounded, and (h1) holds
Then $M=\mathbb{R} \times M_0$ has at least $N$ non-trivial closed geodesics, distinct modulo the action of the group $O(2)$.
Furthermore, if
\begin{itemize}
\item[(h2)] the matrix $[h_{ij}(p,\cdot)]$ representing the bilinear form $h$ is  positive definite for $p \to + \infty$, and negative definite for $p \to - \infty$,
\end{itemize}
then $M$ possesses at least $2N$ non--trivial closed geodesics, geometrically distinct.
\end{teo}

\begin{proof}
Observe that $Z=\mathbb{R} \times Z_0$, where  $Z_0=\{z_{p,q} \mid |p| = |q| =1, \; p \bullet q = 0\}$.
According to Theorem 2.1, it suffices to look for critical points of $\Phi_\ge$. From Lemma \ref{lem:due}, it follows that either $\Phi_\ge = b$ everywhere, or has a critical point $(\bar{r},\bar{p},\bar{q})$. In any case such a critical point gives rise to a (non--trivial) closed geodesic of $(M,g_\ge)$.

From Remark \ref{rem:3}, we know that $\Phi_\ge$ is $O(2)$--invariant. This allows us to introduce the $O(2)$--category $\operatorname{cat}_{O(2)}$. One has
\begin{equation*}
\operatorname{cat}_{O(2)} (Z) \geq \operatorname{cat}(Z/O(2)) \geq \operatorname{cuplength} (Z/O(2)) +1.
\end{equation*}
Since $\operatorname{cuplength}(Z/O(2)) \geq N-1$, (see \cite{sch}), then $\operatorname{cat}_{O(2)} (Z) \geq N$. Finally, by the Lusternik--Schnirel'man theory, $M$ carries at least $N$ closed geodesics, distinct modulo the action $O(2)$. This proves the first statement.

Next, let
\begin{equation}
\Gamma (r,p,q)=G((r,z_{p,q}))=\frac{1}{2}\int_0^1 h(r,z_{p,q}(t))[{\dot z}_{p,q},{\dot z}_{p,q}] \, dt
\end{equation}
Then (h) immediately implies that
\begin{equation}
\Gamma (r,p,q) \to 0 \mbox{\quad as } |r| \to \infty,
\end{equation}
Moreover, if (h2) holds, then $\Gamma (r,p,q) > 0$ for $r > r_0$, and $\Gamma (r,p,q) < 0$ for $r < -r_0$. Since (recall equation \eqref{eq:gamma})
\begin{equation} \label{phi:gamma}
\Phi_\ge (r,p,q)=b+\ge \Gamma (r,p,q) + o(\ge),
\end{equation}
it follows that
\begin{equation*}
\begin{cases}
\Phi_\ge (r,p,q) > b &\text{ for } r > r_0 \\
\Phi_\ge (r,p,q) < b &\text{ for } r < - r_0.
\end{cases}
\end{equation*}
We can now exploit again the $O(2)$ invariance.

By assumption, and a simple continuity argument, $\{ \Phi_\ge > b \} \supset [R_0 , \infty ) \times Z_0$, and similarly $\{ \Phi_\ge < b \} \supset [-\infty , -R_0 ) \times Z_0$, for a suitably large $R_0 > 0$. Hence $\operatorname{cat}_{O(2)} ( \{ \Phi_\ge > b \}) \geq \operatorname{cat}_{O(2)}(Z_0) = N$. The same argument applies to $\{\Phi_\ge < b\}$. This proves the existence of at least $2N$ closed geodesics.
\end{proof}

\begin{rem}
\begin{itemize}
\item[(i)] In \cite{carminati}, the existence of $N$ closed geodesics on $S^N$ endowed with a metric close to the standard one is proved. Such a result does not need any study of $\Phi_\ge$ and its behavior. The existence of $2N$ geodesics is, as far as we know, new. We emphasize that it strongly depends on the {\em form} of $M=\mathbb{R} \times M_0$.

\item[(ii)] In \cite{tanaka}, the metric $g$ on $M$ is possibly not perturbative. No multiplicity result is given.
\end{itemize}
\end{rem}

\section{The general case}

In this section we consider a compact riemannian manifold $(M_0,g_0)$,
and in analogy to the previous section, we put
\begin{equation}
g_\ge (s,\xi)((r,\eta),(r,\eta))=|r|^2+g_0(\xi)(\eta ,\eta)+\ge h(s,\xi)((r,\eta),(r,\eta)).
\end{equation}
Again, we define $\Lambda = \{u=(r,x) \mid r \in H^1(S^1,\mathbb{R}), \; x \in H^1(S^1,M_0) \}$,
\[
E_{M_0}(x)=\frac{1}{2} \int_0^1 g_0(x)(\dot x , \dot x) \, dt,\quad E_0(r,x)=\frac{1}{2} \int_0^1 | \dot r|^2 \, dt + E_{M_0}(x),
\]
and finally
\[
E_\ge (r,x)=E_0(r,x)+\ge G(r,x),
\]
with $G$ as in \eqref{eq:G}.
It is well known (\cite{kli2}) that $M_0$ has a closed geodesic $z_0$. The functional $E_{M_0}$ has again a critical manifold $Z$ given by
\begin{equation*}
Z = \{ u(\cdot)=(\rho,z_0(\cdot + \tau)) \mid \rho \mbox{\; constant, } \tau \in S^1 \}.
\end{equation*}
Let $Z_0 = \{ z_0 ( \cdot + \tau ) \mid \tau \in S^1 \}$. It follows that $Z \approx \mathbb{R} \times Z_0$. The counterpart of $\Gamma$ in \eqref{eq:gamma} is
\begin{equation}\label{eq:gammagenerale}
\Gamma (r,\tau)=\frac{1}{2} \int_0^1 h(r,z_\tau)[\dot z_\tau, \dot z_\tau] \, dt.
\end{equation}

Let us recall some facts from \cite{kli1}.

\begin{rem}
There is a linear operator $A_z \colon T_z \Lambda (M_0) \to T_z \Lambda$, which is a compact perturbation of the identity, such that
\[
E_{M_0}^{\prime \prime} (z) [h,k]=\langle A_z h \mid k \rangle_1 = \int_0^1 \overbrace{A_z h}^{\cdot} \bullet \dot k \, dt.
\]
In particular, $E_0$ satisfies (AS2).
\end{rem}

\begin{defi} \label{defi:nondegen}
Let
\[
\ker E_{M_0}^{\prime \prime}(z_0)=\{ h \in T_{z_0} \Lambda (M_0) \mid \langle A_{z_0} h \mid k \rangle_1 =0 \quad \forall k \in T_{z_0} \Lambda (M_0) \}.
\]
We say that a closed geodesic $z_0$ of $M_0$ is {\em non--degenerate}, if
\[
\dim \ker E_{M_0}^{\prime \prime}(z_0)=1.
\]
\end{defi}

\begin{rem}
For example, it is known that when $M_0$ has negative sectional curvature, then all the geodesics of $M_0$ are non--degenerate. See \cite{DNF}.
Moreover, it is easy to see that the existence of non--degenerate closed geodesics is a {\em generic} property.
\end{rem}

\begin{lem} \label{lem:3.1}
If $z_0$ is a non--degenerate closed geodesic of $M_0$, then $Z$ satisfies (AS2).
\end{lem}

\begin{proof}
It is always true that $T_{z_r} Z \subset \ker E_0^{\prime \prime}(z_r)$. By \eqref{eq:nucleoE0},  we have that  $\dim T_{z_r} Z = \dim \ker E_0^{\prime \prime}(z_r)$. This implies that $T_{z_r} Z = \ker E_0^{\prime \prime}(z_r)$. 
A generic element of $Z$ has the form $(\rho,z^\tau)$ for $\rho \in \mathbb{R}$ and $z^\tau = z(\cdot + \tau)$; then
\[
T_{(\rho ,z^\tau)} M=\mathbb{R} \times T_{z^\tau} M_0,
\]
and any two vector fields $Y$ and $W$ along a curve on $M=\mathbb{R} \times M_0$ can be decomposed into
\begin{equation} \label{eq:Y}
Y=h(t)\frac{d}{dt} + y(t) \in \mathbb{R} \oplus T_{z^\tau} Z_0,
\end{equation}
\begin{equation} \label{eq:W}
W=k(t)\frac{d}{dt} + w(t) \in \mathbb{R} \oplus T_{z^\tau} Z_0.
\end{equation}

In addition, there results (see \cite{hebey})
\begin{equation}\label{eq:14}
E_{M_0}^{\prime \prime} (z_0)[y,w]=\int_0^1 \left[ g_0 (D_t y,D_t w)-g_0(R_{M_0}y(t),\dot z_0 (t))\dot z_0(t) \mid w(t)) \right] \, dt,
\end{equation}
and
\begin{equation}
R_M(r,z)=R_{\mathbb{R}} (r) + R_{M_0} (z) = R_{M_0}(z),
\end{equation}
where $R_M$, $R_{M_0}$, etc.  stand for the curvature tensors of $M$, $M_0$, etc. By \eqref{eq:14}, \eqref{eq:Y} and \eqref{eq:W}, as in the previous section, $E_0^{\prime \prime}(\rho ,z_\tau)[Y,W]=0$ is equivalent to the system
\begin{equation}\label{eq:nondegen}
\begin{cases}
\ddot h =0 \\
\int_0^1 g_0 (z)[D_t y,D_t w]- \langle R_{M_0}(y(t) , \dot z_r(t))\dot z_r(t) \mid w(t) \rangle  \, dt =0.
\end{cases}
\end{equation}
As in the case of the sphere, the first equation implies that $h$ is constant. The second equation in \eqref{eq:nondegen} implies that $y \in \ker E_{M_0}^{\prime \prime} (z^\tau)=\ker E_{M_0}^{\prime \prime}(z_0)$.  Hence,
\begin{equation} \label{eq:nucleoE0}
\ker E_0^{\prime \prime}(z_r)=\{ (h,y) \mid h \mbox{\; is constant, and } y \in \ker E_{M_0}^{\prime \prime} (z_0) \}.
\end{equation}
This completes the proof.
\end{proof}

\vspace{7pt}

\begin{teo} \label{th:3.1}
Let $M \sb 0$ be a compact, connected manifold of dimension $N < \infty$. Assume that $M_0$ admits a non--degenerate closed geodesic $z$, and that (in local coordinates) $h_{ij}(p,\cdot) \to a_{-}$ as $p \to - \infty$, and $h_{ij}(p,\cdot) \to a_{+}$ as $p \to + \infty$.
\begin{enumerate}
\item If $a_{-}=a_{+}$ and $h_{ij}(p, \cdot)$ satisfies (h2), then $M$ has at least one closed geodesic.
\item  If $a_{-} \leq a_{+}$ and $h_{ij}(p,\cdot)[u,v]-a(u \mid v )$ is negative definite for $p \to - \infty$ and positive definite for $p \to + \infty$, then $M$ has at least two non-trivial closed geodesic.
\end{enumerate}
\end{teo}

\begin{proof}
Lemma \ref{lem:3.1} allows us to repeat all the argument in Theorem \ref{th:2.1}, and the result follows immediately.
\end{proof}

\section{Isolated geodesics}

In this final section, we discuss one situaion where the critical manifold $Z$ may be degenerate. 
Here, the non--degeneracy condition (AS3) fails, and $T_z Z \subset \ker E_0^{\prime \prime}(z)$ strictly. Fix a closed geodesic $Z_0$ for $M_0$, and put $\tilde{W}=(T_{z_0} Z)^\bot$. Since $T_z Z \subset \ker E_0^{\prime \prime}(z)$ strictly, there exists $k > 0$ such that $\tilde{W}=(\ker E_0^{\prime \prime}(z_0))^\bot \oplus \mathbb{R}^k$. Repeating the preceding finite dimensional reduction, one can find again a unique map $\tilde{w}=\tilde{w}(z,\zeta)$, where $z \in Z$ and $\zeta \in \mathbb{R}^k$, in such a way that $E_\ge^\prime =0$ reduces to an equation like
\[
\nabla A(z+\zeta + \tilde{w}(z,\zeta))=0.
\]
If $z_0$ is an isolated minimum of the energy $E_{M_0}$ over some connected component of $\Lambda (M_0)$, then it is possible to show that there exists again a function $\Gamma \colon Z \to \mathbb{R}$ such that
\[
\nabla A(z+\zeta + \tilde{w}(z,\zeta))=0 \iff \frac{\de\Gamma}{\de r}(-R,\tau)\frac{\de \Gamma}{\de r}(R,\tau) \ne 0
\]
for some $R \in \mathbb{R}$ and all $\tau \in S^1$. For more details, see \cite{berti}. In particular, we will use the following result.

\begin{teo} \label{th:berti}
Let $H$ be a real Hilbert space, $f_\ge \colon H \to \mathbb{R}$ is a family of $C^2$--functionals of the form $f_\ge = f_0 + \ge G$, and that:
\begin{itemize}
\item[(f0)] $f_0$ has a finite dimensional manifold $Z$ of critical points, each of them being a minimum of $f_0$;
\item[(f1)] for all $z \in Z$, $f_0^{\prime \prime}(z)$ is a compact perturbation of the identity.
\end{itemize}
Fix $z_0 \in Z$, put $W=(T_{z_0}Z)^\bot$, and suppose that $(f_0)_{|W}$ has an isolated minimum at $z_0$. Then, for $\ge$ sufficiently small, $f_\ge$ has a critical point, provided $\deg (\Gamma^\prime , B_R , 0 ) \ne 0$.
\end{teo}

\begin{rem}
Theorem \ref{th:berti} has been presented in a linear setting. For Riemannian manifold, we can either 
reduce to a {\em local} situation and then apply the exponential map, or directly resort to the slightly 
more general degree theory on Banach manifold developed in \cite{ET}.
\end{rem}

\begin{teo}\label{th:4.1}
Assume that $\pi_1 (M_0) \ne \{0\}$, and that all the critical points of $E_0$, the energy functional of $M_0$, are isolated. Suppose the bilinear form $h$ satisfies (h1), and
\begin{itemize}
\item[(h3)] $\dfrac{\de h}{\de r}(R,\xi)\dfrac{\de h}{\de r}(-R,\xi) \ne 0$ for some $R>0$ and all $\xi \in S^1$.
\end{itemize}
Then, for $\ge >0$ sufficiently small, the manifold $M=\mathbb{R} \times M_0$ carries at least one closed geodesic.
\end{teo}

\begin{proof}
We wish to use Theorem \ref{th:berti}. Since $\pi_1 (M_0) \ne \{0\}$, then $E_0$ has a geodesic $z_0$ such that 
$E_0 (z_0) = \min E_0$ over some component $C$ of $\Lambda(M_0)$. See \cite{kli2}.

We consider the manifold
\[
Z = \{ u \in \Lambda \mid u(t)=(\rho,z_0(t+\tau)), \; \rho \mbox{\; constant, } \tau \in S^1 \}.
\]
Here we do not know, a priori, if $Z$ is non--degenerate in the sense of condition (AS2).
But of course $(E_0)_{W}$ has a minimum at the point $(\rho , z_0)$, where
$W=(T_{\rho , z^\tau)}Z)^\bot$. We now check that it 
is isolated for $(E_0)_{W}$. We still know that $Z = \mathbb{R} \times Z_0$. Take any point $(\rho , z_\tau) \in Z$, and observe that 
$T_{(\rho,z^\tau)} Z = \{ (r,y) \mid r \in \mathbb{R}, \; y \in T_{z^\tau}
Z_0 \}$. 
For all $(r,y)\in W$ sufficiently close to $(\rho,z_0)$, it holds in particular that $y \bot z_\tau$. Hence
\[
E_0(r,y) = L_0(r)+E_{M_0}(y) \geq E_{M_0}(y) >
E_{M_0}(z^\tau)=E_{M_0}(z_0) = E_0(\rho , z_0)
\]
since $L_0 \geq 0$ and $z_0$ (and hence $z^\tau$, due to $O(2)$
invariance) is an isolated minimum of $E_{M_0}$ by assumption.

Finally, thanks to assumption (h3), $\frac{\de \Gamma}{\de r}(-R,\tau)\frac{\de\Gamma}{\de r}(R,\tau) \ne 0$.

This concludes the proof.
\end{proof}

\bibliography{geo1}
\bibliographystyle{amsplain}
\end{document}